\documentclass[10pt,a4paper]{article}
\usepackage{comment}
\usepackage{empheq}
\usepackage[T1]{fontenc}
\usepackage[mathscr]{euscript}%pour la police mathscr
\usepackage{amsmath}
\usepackage{mathtools}%pour ecrire au-dessus des rightarrow
\usepackage{extarrows}
\usepackage{amsbsy}
\usepackage{amsfonts}
\usepackage{amssymb}%contient \nexists par exemple
\usepackage{graphicx}
\usepackage{bbm}%pour les indicatrices 1 à double barre
\usepackage{listings}%pour inserer des codes source
\usepackage{color}%fonctionne de pair avec listings
\usepackage[table]{xcolor}%pour les tableaux tabular en couleur
\usepackage{xcolor}
\usepackage[left=1cm,right=2cm,top=2cm,bottom=2cm]{geometry}
\usepackage[latin1]{inputenc}
\usepackage{amsmath,amssymb,amsfonts}
\usepackage[english]{babel}
\usepackage{enumerate}
\usepackage{yhmath}%pour ring et widering
\usepackage{amsthm}
\usepackage{picture}
\usepackage{colortbl}%pourcolorier les cellules des tableaux
\usepackage{hyperref}
\usepackage{wrapfig} %pour qu'une figure soit encadree par le texte
\usepackage{graphicx}%pour inclure des images dans le texte
\usepackage{xfrac}  %pour les quotients \sfrac
\usepackage{faktor}  %pour les quotients \faktor
\usepackage{pst-all} %pour utiliser PSTricks
\usepackage{pstricks}
\usepackage{permute} % pour les permutations pmt

\usepackage{array}
\usepackage{tikz}
\usepackage{comment} 
\usepackage{dsfont} %pour les lettres minuscules en mathbb
\usepackage{amsmath}
\usepackage{amsfonts}
\usepackage{amssymb}

\definecolor{ashgrey}{rgb}{0.7, 0.75, 0.71}%cf latexcolor.com

\usetikzlibrary{arrows}
\usetikzlibrary{arrows,shapes,positioning}
\usetikzlibrary[patterns]
\usetikzlibrary{decorations.markings}
\usepackage{pgf,tikz}
\usepackage{pgfplots}
\tikzstyle directed=[postaction={decorate,decoration={markings,
    mark=at position .5 with {\arrow[]{stealth}}}}]
\usepackage{tikz-cd} %pour les diagrammes commutatifs
\usetikzlibrary{cd}
\usetikzlibrary{calc}
\usetikzlibrary{trees,positioning,arrows,fit,shapes}
\usetikzlibrary{matrix} %pour les diagrammes commutatifs et matrices et les fonctions f:A->B
\usetikzlibrary{angles,quotes,calc}
\definecolor{gris-clair}{rgb}{0.8,0.8,0.8}

\tikzstyle directed=[postaction={decorate,decoration={markings, 
mark=at position .65 with {\arrow{latex}}}}]%pour que les pointes des fleches soient au milieu

\usepackage{tikz-cd}

\tikzset{
  symbol/.style={
    draw=none,
    every to/.append style={
      edge node={node [sloped, allow upside down, auto=false]{$#1$}}}
  }
}%pour des fleches en symbole (voir le chapitre sur l'algebre de decomposition d'un polynome unitaire)

\usepackage[all]{xy}
\def\dar[#1]{\ar@<1pt>[#1]\ar@<-1pt>[#1]}
\entrymodifiers={!!<0pt,0.7ex>+}

\usepackage{amsthm} %pour les theoremes, lemmes et preuves
\theoremstyle{definition}

\theoremstyle{definition}
\title{Solutions in $\mathbb{Z}[i]$ of $A^{5}+B^{5}=C^{5}\pm 1$}
\author{D. FOSSE, MSc. Physics\\
dominique.fosse@a3.epfl.ch}\date{\vspace{-5ex}}
\begin{document}
\maketitle
\paragraph{Introduction}
In \cite{Hirschhorn}, it is shown how we can derive infinitely many integral solutions of $A^{3}+B^{3}=C^{3}\pm 1$ starting from Ramanujan's identity:
\begin{equation*}
{(a^{2}+7ab-9 b^{2})}^{3}+{(2a^{2}-4 ab+12 b^{2})}^{3}={(2 a^{2}+10 b^{2})}^{3}+{(a^{2}-9 ab -b^{2})}^{3}
\end{equation*}
We will do something similar but for the equation $A^{5}+B^{5}=C^{5}\pm 1$; the price to pay will be that the proposed solutions are Gaussian integers and not natural integers. With $i^{2}=-1$, we consider: 
\begin{equation*}
g(x):=(x^{2}+2ax-2 a^{2})^{5}+(i x^{2}-2ax+2i a^{2})^{5}
\end{equation*}
One gets directly $g(x)=g(-x)$ so $g$ is even and this means:
\begin{equation}\label{eq-1}
(x^{2}+2ax-2 a^{2})^{5}+(i x^{2}-2ax+2i a^{2})^{5}=(x^{2}-2ax-2 a^{2})^{5}+(i x^{2}+2ax+2i a^{2})^{5}
\end{equation}
With $(a,x)\in\mathbb{Z}^{2}$, one gets Gaussian integers solutions of $A^{5}+B^{5}=C^{5}+D^{5}$ which is a good starting point.
\paragraph{Using classical linear recurrence relations} Let's consider the following linear recurrence, $\forall n\in\mathbb{N}$:
\begin{equation}\label{eq-2}
F_{n+2}=-2 F_{n+1}+2 F_{n}\qquad F_{0}:=0,\, F_{1}:=1
\end{equation}
This can be solved using classical linear algebra:
\begin{equation}\label{eq-3}
F_{n}=\frac{\sqrt{3}}{6}\left(\left(-1+\sqrt{3}\right)^{n}-\left(-1-\sqrt{3}\right)^{n}\right)\qquad\forall n\in\mathbb{N}
\end{equation}
From this, we compute by direct calculation:
\begin{align*}
F_{n+1}^{2}&=\frac{2^{n}}{6}\left(\left(2+\sqrt{3}\right)^{n+1}+\left(2-\sqrt{3}\right)^{n+1}+2 (-1)^{n} \right)\\
F_{n+1}F_{n}&=\frac{1}{12}\left(\left(-1+\sqrt{3}\right)^{2n+1}+\left(-1-\sqrt{3}\right)^{2n+1}-(-2)^{n+1} \right)\\
F_{n+2}F_{n}&=\frac{2^{n}}{6}\left(\left(2+\sqrt{3}\right)^{n+1}+\left(2-\sqrt{3}\right)^{n+1}+4 (-1)^{n+1} \right)
\end{align*}
From this, we get, in particular, 
\begin{equation}\label{eq-4}
F_{n+1}^{2}-F_{n}F_{n+2}=2^{n}(-1)^{n}
\end{equation}
Let's take $F_{n+1}:=x$ and $F_{n}:=a$ and replace this in equation \eqref{eq-1}, one gets:
\begin{align}
A_{n}:=& F_{n+1}^{2}-2 F_{n+1}F_{n}-2 F_{n}^{2}\nonumber\\
=&\frac{2^{n}}{3}\left((1+\sqrt{3})(2+\sqrt{3})^{n}+(1-\sqrt{3})(2-\sqrt{3})^{n}+(-1)^{n}\right)\\
B_{n}:=& i F_{n+1}^{2}+2 F_{n+1}F_{n}+2i F_{n}^{2}\nonumber\\
=&\frac{2^{n}}{6}\left((-1+\sqrt{3})(2-\sqrt{3})^{n}-(1+\sqrt{3})(2+\sqrt{3})^{n}+2(-1)^{n}\right)+i\frac{2^{n}}{6}\left((3-\sqrt{3})(2-\sqrt{3})^{n}+(3+\sqrt{3})(2+\sqrt{3})^{n}\right)\\
C_{n}:=& i F_{n+1}^{2}-2 F_{n+1}F_{n}+2i F_{n}^{2}\nonumber\\
=&\frac{2^{n}}{6}\left((1+\sqrt{3})(2+\sqrt{3})^{n}+(1-\sqrt{3})(2-\sqrt{3})^{n}-2(-1)^{n}\right)+i\frac{2^{n}}{6}\left((3-\sqrt{3})(2-\sqrt{3})^{n}+(3+\sqrt{3})(2+\sqrt{3})^{n}\right)\\
d_{n}:=& F_{n+1}^{2}+2 F_{n+1}F_{n}-2 F_{n}^{2}\nonumber
\end{align}
Thus, equation \eqref{eq-1} becomes $A_{n}^{5}+B_{n}^{5}=C_{n}^{5}+d_{n}^{5}$. Now, 
\begin{align*}
d_{n}=F_{n+1}^{2}+2 F_{n+1}F_{n}-2 F_{n}^{2}&=F_{n+1}^{2}- F_{n}\left(-2F_{n+1}+2F_{n}\right)\\
&=F_{n+1}^{2}- F_{n}F_{n+2}\qquad\textnormal{using equation \eqref{eq-2}}\\
&=2^{n}(-1)^{n}\qquad\textnormal{using equation \eqref{eq-4}}
\end{align*}
With $z\in\mathbb{C}$, $|z|<1$, it's easy enough to calculate the following, from equations $(5)$, $(6)$ and $(7)$:
\begin{align*}
\sum_{n\geq 0}A_{n}z^{n}&=\frac{4 z^{2}+1}{(2z+1)(4 z^{2}-8z+1)}\\
\sum_{n\geq 0}B_{n}z^{n}&=\frac{-4z}{(2z+1)(4 z^{2}-8z+1)}+i\left(\frac{1-2 z}{4 z^{2}-8z+1}\right)\\
\sum_{n\geq 0}C_{n}z^{n}&=\frac{4z}{(2z+1)(4 z^{2}-8z+1)}+i\left(\frac{1-2 z}{4 z^{2}-8z+1}\right)
\end{align*}
And this satisfies $A_{n}^{5}+B_{n}^{5}=C_{n}^{5}+2^{5n}(-1)^{n}\implies\left(\frac{A_{n}}{2^{n}}\right)^{5}+\left(\frac{B_{n}}{2^{n}}\right)^{5}=\left(\frac{C_{n}}{2^{n}}\right)^{5}+(-1)^{n}$. Let's take $x:=2z$ and $\forall n\in\mathbb{N}$, $a_{n}:=\frac{A_{n}}{2^{n}}$, $b_{n}:=\frac{B_{n}}{2^{n}}$ and $c_{n}:=\frac{C_{n}}{2^{n}}$. It follows that
\begin{empheq}[box=\fbox]{align*}\label{eq-8}
\sum_{n\geq 0}a_{n}x^{n}&=\frac{x^{2}+1}{(x+1)(x^{2}-4x+1)}\\
\sum_{n\geq 0}b_{n}x^{n}&=\frac{-2x}{(x+1)(x^{2}-4x+1)}+i\left(\frac{1-x}{x^{2}-4x+1}\right)\\
\sum_{n\geq 0}c_{n}x^{n}&=\frac{2x}{(x+1)(x^{2}-4x+1)}+i\left(\frac{1-x}{x^{2}-4x+1}\right)\\
\implies& a_{n}^{5}+b_{n}^{5}=c_{n}^{5}+(-1)^{n}
\end{empheq}
For example, 
\begin{align*}
n:=1\implies& 3^{5}+(-2+3 i)^{5}=(2+3 i)^5-1\\
n:=2\implies& 13^{5}+(-6+11 i)^{5}=(6+11 i)^{5}+1\\
n:=3\implies& 47^{5}+(-24+41 i)^{5}=(24+41 i)^{5}-1
\end{align*}

\end{document}